\documentclass[final,3p,times]{elsarticle}
\usepackage{graphicx}
\usepackage{subfigure}
\usepackage{amssymb}
\usepackage{amsfonts}
 \usepackage{lineno}
\usepackage{amsmath}
\usepackage{mathrsfs}
\usepackage{dcolumn}
\usepackage{comment}
\usepackage{color}
\usepackage{float}

\newcommand{\beq}{\begin{equation}}
\newcommand{\eeq}{  \end{equation}}
\newcommand{\beqa}{\begin{eqnarray}}
\newcommand{\eeqa}{  \end{eqnarray}}

\graphicspath{{../Figs/}}

\def\XXint#1#2#3{{\setbox0=\hbox{$#1{#2#3}{\int}$}
     \vcenter{\hbox{$#2#3$}}\kern-.5\wd0}}

\definecolor{red}{rgb}{1,0,0}
\definecolor{blue}{rgb}{0,0,1}

\def\strutdepth{\dp\strutbox}
\def\nw#1{\strut\vadjust{\kern-\strutdepth\vtop to0pt{\vss\hbox to\hsize
{\hskip\hsize\hskip5pt$\leftarrow$\hss\strut}}}{\em #1}}

\begin{document}

\begin{frontmatter}
	\title{MARS~: a Method for the Adaptive Removal of Stiffness in PDEs
}


\author[LD]{Laurent Duchemin}
\ead{laurent.duchemin@espci.fr}
\author[JE]{Jens Eggers}
\ead{Jens.Eggers@bristol.ac.uk}

	\address[LD]{Physique et M\'ecanique des Milieux H\'et\'erogenes, CNRS, ESPCI Paris, Universit\'e PSL, Sorbonne Universit\'e, Universit\'e de Paris, F-75005 Paris, France}
\address[JE]{Department of Mathematics - University of Bristol, University Walk, Bristol BS8 1TW, United Kingdom}


\begin{abstract}
The E(xplicit)I(implicit)N(null) method was developed recently to remove numerical instability
from PDEs, adding and subtracting an operator $\mathcal{D}$ of arbitrary 
structure, treating the operator implicitly in one case, and explicitly in the other. 
Here we extend this idea by devising an adaptive procedure to find an 
optimal approximation for $\mathcal{D}$. We propose a measure of the numerical
error which detects numerical instabilities across all wavelengths, and adjust
each Fourier component of $\mathcal{D}$ to the smallest value such that 
numerical instability is suppressed. We show that for a number of nonlinear
and non-local PDEs, in one and two dimensions, the spectrum 
of $\mathcal{D}$ adapts automatically and dynamically to the theoretical 
result for marginal stability. Our method thus has the same stability properties as 
a fully implicit method, while only requiring the computational cost comparable to an explicit solver. 
The adaptive implicit part is diagonal in Fourier space, and thus leads to 
minimal overhead compared to the explicit method. 
\end{abstract}


\end{frontmatter}




\section{Introduction}

Our ability to model many key physical processes is limited 
by the stability of the numerical schemes we use to simulate the
partial differential equations (PDEs) describing them.
The reason is that the maximal stable time step of an explicit 
numerical integration scheme is of the order of the shortest time-scale
in the system. In a stable physical system these are typically exponentially
damped modes which relax back to equilibrium; the smaller the length
scale, the faster the relaxation. This makes it particularly hard to simulate
systems at large values of the viscosity or of the surface tension. 
For instance, surface tension driven flows in the open source fluid dynamics
code Gerris \cite{Pop09} (followed by Basilisk~: \verb http://basilisk.fr ) 
require a time step proportional to
$\Delta^{3/2}$\cite{Brackbill1992}, where $\Delta$ is the grid spacing,
which this program adapts dynamically in order to ensure a sufficient
spatial accuracy \cite{Popinet2018}. As a result, for small geometries 
$\Delta$ can be very small, resulting in time steps which are prohibitively
small. This constraint is more restrictive than the CFL constraint,
related to advection, for which the time step depends on scale like $\Delta$.
Another example is the numerical computation of
solidification/fusion fronts, which uses a non-linear heat
equation \cite{Labrosse2012}~: the corresponding time step constraint
is $\Delta^2$.

If for example relaxation toward equilibrium is controlled
by a differential operator of order $m$ ($m=2$ for ordinary diffusion, 
$m=3$ for the Hele-Shaw flow to be described below),
then the required maximum time step $\delta t$ scales as 
$\delta t = C \delta x^m$, where $\delta x$ is the smallest grid spacing
or the size of the smallest sub-division. In a well-resolved numerical 
simulation, this should be considerably smaller than the smallest relevant 
physical feature. Rapid exponential decay implies that the amplitude of 
perturbations on the grid scale is very small, and contributes negligible 
to the numerical solution. Thus one arrives at the paradoxical situation 
that the stability of the numerical scheme is controlled by a part of the 
solution which contributes negligibly, and which is actually the most
stable from a physical perspective. This property is sometimes referred
to as the stiffness of the PDE \cite{KT05}, which becomes worse with increasing
spatial order $m$ of the operator. 

To deal with this constraint on the time step, which often is so severe
that it makes the exploration of important physical parameter regimes
impractical, one has to resort to implicit
methods. This means that the right hand side of the equation (or at least
the stiffest parts of it) has to be evaluated at a {\it future} time step,
making it necessary to solve an implicit equation at each time
step \cite{Iserles_book,Ames_book}. This makes the numerical code both
complicated to write and time-consuming to solve. 
This is true in
particular if the operator is non-local (as is the case for example
of integral operators, as they appear in boundary integral type
codes \cite{Pozbook,HLS94}). Indeed, in this particular case, when
writing an implicit scheme, 
each element of the discretized solution depends on all the others, 
requiring a large number of operations to solve the implicit equation. 

To address this problem, it has long been realized that not the whole of
the right hand side of an equation has to be treated implicitly, as long
as the ``stiffest'' part of the operator is dealt with implicitly. This
gives rise to the so-called ``implicit-explicit methods'' \cite{ARW95,DB12},
which divide up the problem between explicit and implicit parts, such
that hopefully the implicit contribution is sufficiently simple to invert.
If this is not clear, as is typically the case for an integral operator,
the problem can be solved by judiciously slicing off the stiffest part,
which can be local \cite{HLS94}). However, this has to be done on a
case-by-case basis, and will not always be possible. Recently, we have
presented a much more general method to stabilize stiff equations, which
makes use of the arbitrariness in which splitting between explicit and
implicit parts can take place \cite{ELS99,DE14}. We consider a partial
differential equation of the form
\beq
\frac{\partial u}{\partial t} = f(u,t),
\label{edp}
\eeq
where $u(x,t)$ is a function of space and time or a vector of functions
of space and time, and $f(u,t)$ generally is a non-linear operator
involving spatial derivatives of $u(x,t)$. 
In the present article, after explaining the adaptive stabilization procedure in section \ref{sec:dynamic}, 
we shall treat the following three examples~: 
\begin{itemize}
\item A non-linear operator with a fourth-order spatial derivative,
  related to the thin film flow equation (section \ref{sec_film})~:
		$$
		f(u,t) = 
- \frac{\partial}{\partial x} \left( u^3
\frac{\partial^3 u}{\partial x^3} 
+\frac{1}{u}\frac{\partial u}{\partial x}
\right),
		$$
	\item A two-dimensional example with a fourth-order derivative (section \ref{sec_kur})~:
		$$
		f(u,t) = 
  		-\mathcal{N}(u) - \Delta u - \nu \Delta^2 u, 
		$$
where $u(x,y,t)$, $\mathcal{N}(u)$ is a non-linear operator,
$\Delta$ the Laplacian, and $\nu$ a constant, 
	\item A boundary integral equation (section \ref{sec_HS})~:
		$$
		f(u,t) = \int{g(v) K(u,v)} dv,
		$$
		where $u$ and $v$ are a two-dimensional vectors, $g(v)$
                is a function of $v$ involving second-order derivatives
                in space, $K(u,v)$ is a singular kernel, and the integral
                is performed along a curve. 
\end{itemize}

As explained in our previous article \cite{DE14}, in order to stabilize the stiff terms in $f(u,t)$,
we add two terms on the right-hand-side of the discretized version of
equation (\ref{edp})~:
\beq
\frac{u_j^{n+1}-u_j^n}{\delta t} = f_j(u^n,t^n) 
- \mathcal{D}_j[u^n] 
+ \mathcal{D}_j[u^{n+1}],
\label{discrete}
\eeq
where $n$ denotes the time variable ($t^n=n \delta t$) and $\mathcal{D}$
is an arbitrary operator.
The variable $u$ as well as $f$ are defined on a spatial grid
$x_j=j\delta x$, where $\delta x$ is the grid spacing. Clearly, the added
terms are effectively zero apart from the first-order error that comes
from the fact that $\mathcal{D}$ is evaluated at different time levels,
which motivates the name ``Explicit-Implicit-Null'' method or ``EIN''. If
$\mathcal{D}$ is the same as the original operator $f(u,t)$, this is a
purely implicit method, if $\mathcal{D}=0$, it is explicit. 
Similar ideas have been implemented to stabilize the motion of a surface
in the diffuse interface and level-set methods
\cite{smereka02,glasner03,SL08}, and for the solution of PDEs on surfaces
\cite{MR09}. We also show that by a simple step-halving procedure
\cite{AD04}, (\ref{discrete}) can always be turned into a scheme which
is second order accurate in time \cite{DE14}. 

When applying \eqref{discrete}, we want $\mathcal{D}$ to be a reasonable
approximation to the stiff part of $f$. If $\mathcal{D}$ were much larger,
it would stabilize the scheme, but would introduce an additional time
truncation error. We thus want to choose $\mathcal{D}$ for optimal 
effectiveness, in the sense that it adds the perfect amount of damping, 
without adding a supplementary error to the numerical scheme.

This paper presents a numerical scheme to achieve this goal automatically,
adjusting $\mathcal{D}$ to the threshold value. Restricting ourselves to
periodic boundary conditions, we choose $\mathcal{D}$ to be diagonal in
Fourier space, which renders it both simple to handle and sufficiently
flexible. Indeed, the implicit step becomes almost trivial to perform:
\beq
\frac{\hat{u}^{n+1}_k-\hat{u}^n_k}{\delta t} =
\hat{f}_k(u^n,t^n) + \lambda(k) \hat{u}^n_k 
- \lambda(k) \hat{u}^{n+1}_k ,
\label{discrete_F}
\eeq
where ~$\hat{}$~ denotes the Fourier transform and the damping spectrum
$\lambda(k) \ge 0$ is an arbitrary function. Adjustment of $\mathcal{D}$
is now down to finding the scalar damping spectrum $\lambda(k)$ for a
discrete sequence of $k$'s. 

The Fourier transform
$\hat{f}_k$ can be calculated effectively from the spatial discretization
$f_j$ using the fast Fourier transform (FFT) \cite{FFTW05}. From
\eqref{discrete_F}, we find
\beq
\hat{u}^{n+1}_k = \hat{u}^n_k +
\frac{\hat{f}_k(u^n,t^n)}{\delta t^{-1}+\lambda(k)},
\label{discrete_F_solve}
\eeq
so we obtain the desired solution ${u_j}^{n+1}$ at the new time step 
from the inverse transform. The scheme (\ref{discrete_F_solve}) (as well as
any other first order scheme) can be turned into a second order scheme
by Richardson extrapolation \cite{AD04}. Namely, let $u^{1,n+1}$ be the solution
for one step $\delta t$, $u^{2,n+1}$ the solution for two half steps
$\delta t/2$. Then 
\beq
u^{n+1} = 2u^{2,n+1} - u^{1,n+1} + \mathcal{O}(\delta t^3),
\label{R}
\eeq
is second order accurate in time, and 
\beq
E = u^{1,n+1} - u^{2,n+1}
\label{error}
\eeq
can be used as an error estimator \cite{hairer_book}. 

To analyze (\ref{discrete_F}) further, we adopt a ``frozen-coefficient''
hypothesis, that the solution is essentially constant over the time scale
on which numerical instability is developing. Then assuming small perturbations
$\delta \hat{u}_k^{n}$, the problem is turned into a linear equation for
$\delta \hat{u}_k^{n}$, with constant coefficients. At least on a small
scale (i.e. in the large $k$ limit), much smaller than any externally
imposed scale, $\hat{f}_k(u^n,t^n)$ is expected to be translationally
invariant, making the operator diagonal in Fourier space. Namely,
let us assume a more general non-local operator
\[
f(x) = \int_{-\infty}^{\infty} D(x,y) u(y) dy. 
\]
Translational invariance implies that $D(x,y) = D(x-y)$; taking the
Fourier transform, we arrive at $\hat{f}(k) = \hat{D}(k)\hat{u}(k)$. 
This means in the large $k$ limit we expect
\beq
\hat{f}_k(u^n+\delta u,t^n) \sim - e(k) \delta \hat{u}^n_k, 
\label{linearize}
\eeq
where $\delta u$ is a small perturbation around the solution at $t^n$.
Here we assume that the eigenvalues $e(k)$ are real, as it
is typically the case for physical problems, where the dominant process on
a small scale is dissipative. 

We have shown in \cite{DE14} that, as long as $\lambda(k) > e(k)/2$, the
system (\ref{discrete_F}) with the approximation \eqref{linearize}
is unconditionally stable. This is a generalization
of a method first presented, for the case of the diffusion equation in two
dimensions, in \cite{DD71}. If (\ref{discrete_F_solve}) is turned into a
second order scheme using (\ref{R}), this condition is \cite{CG73}~:
\beq
\lambda(k) > \lambda_c(k) = \frac{2}{3}e(k),
\label{2ndstab}
\eeq
with $\lambda_c(k)$ the theoretical stability limit.
Thus for sufficiently large values of $\lambda(k)$, there is always stability;
however, the time truncation error of the second order scheme is
$O(\lambda\delta t^2)$. This is to be compared to the time truncation
error $O(e \delta t^2)$ of the original scheme, and should therefore
be kept as small as possible, consistent with the stability constraint.
There is a certain similarity here with the preconditioning of matrices,
where a matrix is approximated by a simple diagonal matrix
\cite{Wathen1993,ESW_book}.

We also showed in \cite{DE14} that the explicit scheme ({\it i.e.}
for which $\lambda(k)=0$) is stable as long as~:
\beq
e(k) < \frac{2}{\delta t}.
\label{def_ke}
\eeq
Equation \eqref{def_ke} defines a threshold value of wave numbers $k_e$,
below which the scheme is stable even without stabilization. As a result,
for $k < k_e$, $\lambda(k)$ can be chosen to vanish, without affecting
stability. 

In \cite{DE14} we have tested the ideas underlying the EIN method,
calculating the spectrum $e(k)$ for a variety of operators, including
nonlocal operators treated previously in \cite{HLS94}.
We approximated $\lambda(k)$ as a power law, derived from the low
wavenumber limit of the exact discrete spectrum. As predicted by the
above analysis, we find the scheme (\ref{discrete_F}) unconditionally
stable, and performing with the same accuracy as that proposed in
\cite{HLS94}. Obviously, this still requires one to obtain a good
estimate for the spectrum. 

In the present paper, we aim to remove this analytical step, and to
make the calculation of $\lambda(k)$ self-consistent. The idea is to
determine $\lambda(k)$ iteratively, by detecting numerical instability.
If, for a given wave number $k$, random perturbations due to numerical 
instability of the time-stepping grow in time, then the damping is increased,
while $\lambda(k)$ can be reduced if the code is stable. In the simplest
version of our
procedure, we focus on the high wave number limit, where most of the
stiffness is coming from, and approximate $\lambda(k)$ by a power law,
determined by one or two parameters, depending on whether the exponent
is to be prescribed.
While we found this approach to work, it introduces arbitrary assumptions
into the procedure, and assumes a separation between a high and low wave
number regimes. Instead, here we present the results of a scheme which
adjusts each Fourier mode individually, based on noise detected in the
same Fourier mode. This models the original operator in much greater
detail, and leads to a spectrum
$\lambda(k)$ which corresponds closely to the theoretical stability limit. 

In the next section we develop and describe our procedure for
automatic stabilization. The following three sections are each
dedicated to a particular example, to illustrate how the method is
implemented, and to demonstrate its effectiveness. 

\section{Adaptive stabilization}
\label{sec:dynamic}

Our method is based on the formulation (\ref{discrete_F}), which together
with (\ref{R}) is an unconditionally stable second order scheme, as long
as $\lambda(k)$ is sufficiently large. We would like to find an adaptive
procedure which refines $\lambda(k)$ at each time step, so as to keep it
as small as possible, consistent with stability. To achieve this, we have
to address two issues:
(i) find a measure $\epsilon(k)$ of the noise, or of numerical instability,
for each Fourier mode $k$; 
(ii) specify the evolution of $\lambda(k)$ for a given noise. 

Finding a suitable measure of the error is the crucial question,
to be discussed in more detail below. As for (ii), we aim to adjust each
Fourier component $\lambda(k)$ individually, although we have also
explored representing $\lambda(k)$ by a finite number of parameters. 
We adopt a simple approach, taking a local relation between $\epsilon(k)$
and $\lambda(k)$, which is shown to be sufficient for the examples to
be presented below. For each Fourier mode, if $\epsilon(k)$ is larger
than an upper bound $\epsilon_u$, the corresponding $\lambda(k)$ is
increased by a factor of 1.2. If on the other hand $\epsilon(k) < \epsilon_u$,
$\lambda(k)$ is decreased slowly by a factor of $1/1.02$ at each time step, 
in order to avoid a sudden onset of instability. We have used the same
rates in all examples, but confirmed that the method is robust against
change of parameters. 

\begin{figure}[htbp]
\begin{center}
\includegraphics[width=\linewidth]{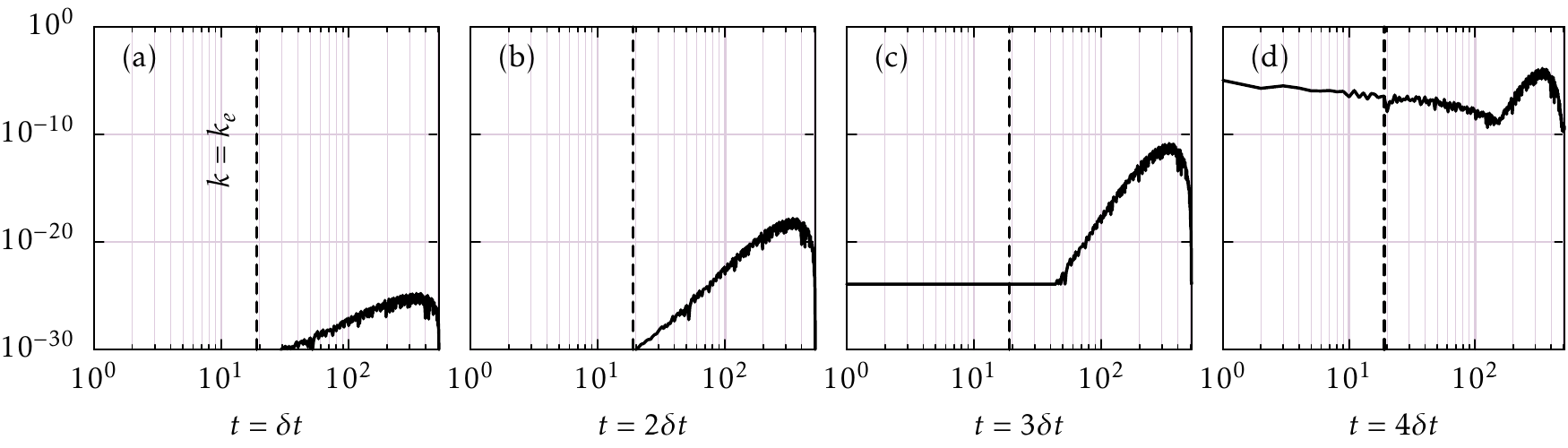}
\caption{
The evolution of the error estimator $\hat{E}_k$ (black curve)
for the Hele-Shaw flow \eqref{advection}, \eqref{BR}, discussed in
more detail in Sect.~\ref{sec_HS}. Initially the error is uniformly small
for a flat interface with a white noise. 
The vertical dashed line is the stability boundary $k=k_e$.
   }
\label{define_epsilon2}
\end{center}
\end{figure}

As to a measure of noise, a first guess might be to take $\epsilon(k)$
as the Fourier transform of the error estimator \eqref{error} $\hat{E}_k$.
We tested this idea using the interface dynamics discussed in more detail in 
section \ref{sec_HS}, and illustrated in Fig.~\ref{HS_stable}.
Figure \ref{define_epsilon2} shows the evolution of the Fourier
transform $\hat{E}_k$ of this error estimator, for the first four
time steps, without using the EIN method.
The time step is chosen to be $\delta t = 3.125\times 10^{-5}$, the number
of points $N=1024$, and we use a purely explicit scheme (no stabilization),
so that the modes with the largest wavenumbers are unstable.

Indeed, as explained in the next sections, there exists a region $k>k_e$ in
$k$-space which is stable with an explicit scheme (on the left of the vertical
dashed line), and an unstable region (on the right), where we would like
to detect numerical instability. As a result, the noise level grows very
rapidly for the right-hand side of the spectrum, and for the first two
time steps there is little power in the $k<k_e$ modes. Thus $\hat{E}_k$
could be used to detect correctly the numerical instability for
large $k$.

However, the left part of the spectrum is soon invaded through non-linear
mode-coupling, and there grows a considerable component of the error at
small $k$ (corresponding to large scales), which would not be damped away
if $\lambda(k)$ was increased. 
The problem is clear: in the proposed scheme, there is no clean distinction
between noise resulting from numerical instability, and the broad spectrum
of unstable modes which is part of the physical solution. The crucial
problem of defining the numerical noise $\epsilon$ lies in this distinction. 

A successful
procedure came from the idea of spatial smoothing, taking the truncation
error as the starting point. To compute $\epsilon$ at the j'th gridpoint,
we consider the error estimator $E_j$, and compare it to a smoothed version
$\bar{E}_j$ at the same point. The reasoning is that $\bar{E}_j$ contains
the full spectrum coming from the deterministic nonlinear dynamics,
so $E_j - \bar{E}_j$ only contains the random noise produced by numerical
instability.  
There are many possible choices for the smoothed-out error.
We chose a polynomial approximation over $2n$ gridpoints, but excluding
$j$ itself, otherwise $\epsilon$ would be identically zero. In other words,
\begin{equation}
\bar{E}_j = {\cal P}\left(E_{j-n},\dots,E_{j-1},E_{j+1},\dots,E_{j+n}\right), 
	\label{interp}
\end{equation}
where ${\cal P}$ is the $(2n-1)$ degree polynomial, passing through
$\left(E_{j-n},\dots,E_{j-1},E_{j+1},\dots,E_{j+n}\right)$. 
Taking the Fourier transform of the difference between $E_j$ and $\bar{E}_j$,
we define the noise measure $\epsilon(k)$ as
\begin{equation}
\epsilon(k) = \hat{E}_k - \hat{\bar{E}}_k.
\label{epsilon}
\end{equation}

\begin{figure}[htbp]
\includegraphics[width=\linewidth]{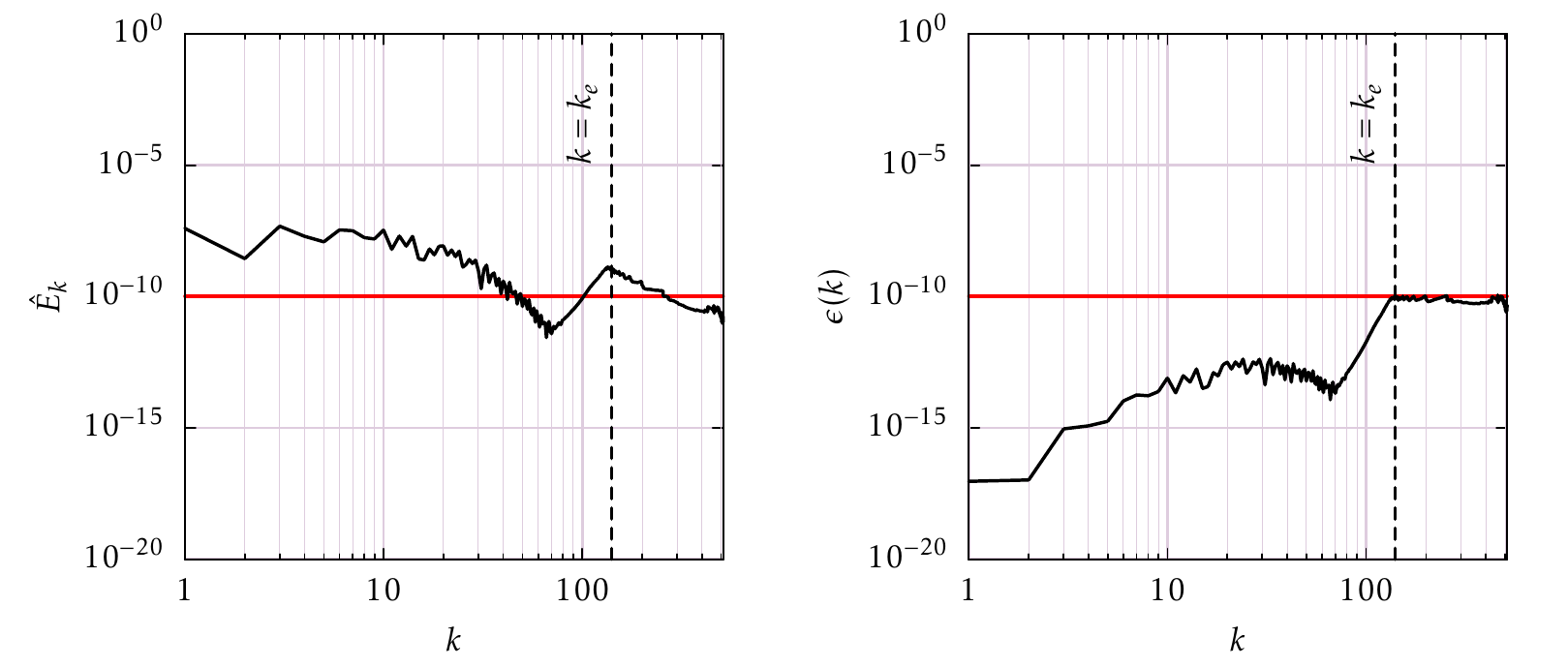}
\begin{center}
\caption{The effect of smoothing on the solution of the Hele-Shaw flow
\eqref{advection},\eqref{BR} at $t=0.04$, as shown in the third panel
of Fig.~\ref{HS_stable} below. On the left, we show the spectrum of the
error estimator $\hat{E}_k$, which is broad in the nonlinear regime. On the 
right, we show the noise measure $\epsilon(k)$ as defined by 
(\ref{epsilon}), which is substantial only in a high wave number region
where noise is detected. The vertical dashed line is $k=k_e$. 
The horizontal red line is the threshold $\epsilon_u$ used to adapt
$\lambda(k)$.
   }
\label{noise_detect}
\end{center}
\end{figure}

The difference between the naive error measure $\hat{E}_k$ and $\epsilon(k)$
based on smoothing is illustrated in Fig.~\ref{noise_detect}. 
We ran the same computation as in figure 
\ref{define_epsilon2}, but using our adaptive procedure. 
Equation \eqref{epsilon} is used as a measure of the noise (with $n=2$)
to adapt $\lambda(k)$ at each time step, with the threshold
$\epsilon_u=10^{-10}$. 
The left curve shows the Fourier transform $\hat{E}_k$ as a function of
$k$~: clearly, this error alone is ill-suited to detect instability,
since it has significant components for $k<k_e$, where the explicit scheme
is stable, {\it i.e.} where there is no instability even for $\lambda(k)=0$. 
The right curve shows the noise measure $\epsilon(k)$ given by
\eqref{epsilon} used to adapt
$\lambda(k)$ as a function of $k$~: the instability is correctly detected
at large values of $k$ and this error remains low for $k<k_e$, {\it i.e.}
does not require any damping, in the region where an explicit scheme is stable.


\section{Example: thin film flow with van der Waals forces}
\label{sec_film}
\subsection{Equation of motion}

As an example of a non-linear equation in one dimension, we first consider
a thin liquid film on a horizontal solid substrate. Assuming lubrication
theory and taking into account van der Waals forces, which can destabilize
the film, the 1D evolution equation for the height of the film $h(x,t)$
reads \cite{williams82,lister99}~:
\begin{equation}
\frac{\partial h}{\partial t} = 
- \frac{\partial}{\partial x} \left( \frac{h^3}{3\eta}
\frac{\partial}{\partial x}
\left[
\gamma \frac{\partial^2 h}{\partial x^2} 
-\frac{A}{6\pi h^3}
\right]
\right), 
\label{thin_film_dim}
\end{equation}
where $\eta$ is the dynamic viscosity of the fluid, $\gamma$ its
surface tension coefficient and $A$ the Hamaker constant. 
Using $L=\sqrt{A/2\pi \gamma}$ as the lengthscale and
$T=3\eta L / \gamma$ as the timescale, the dimensionless equation reads~:
\begin{equation}
\frac{\partial h}{\partial t} = 
- \frac{\partial}{\partial x} 
\left( 
h^3
\frac{\partial^3 h}{\partial x^3} 
+\frac{1}{h}\frac{\partial h}{\partial x}
\right).
\label{thin_film}
\end{equation}

Considering the linear stability of \eqref{thin_film}, we study the
growth of a small-amplitude single mode added to an initially flat interface~:
\begin{equation}
h(x,t) = h_0 + \varepsilon e^{i k x + \omega t},
\label{h_lin}
\end{equation}
where $x\in[0,1]$. Linearizing equation \eqref{thin_film} for
$\varepsilon \ll 1$ gives the dispersion relation~:
\begin{equation}
	\omega = -h_0^3 k^4 + \frac{k^2}{h_0}.
	\label{}
\end{equation}

As the initial condition, we start from a flat film with a sinusoidal
perturbation added to it, corresponding to the most unstable (or Rayleigh)
mode. In order to initiate an instability on the scale of the entire
computational domain, we fix $k=2\pi$ and set the initial height $h_0$
that corresponds to the maximum growth rate~:
$$
\left. \frac{d\omega}{dk}\right|_{k=2\pi}=0 
\iff h_0 = \frac{1}{2^{1/4}(2\pi)^{1/2}} \simeq 0.34
$$
Using this initial thickness, we choose as the initial condition~:
\begin{equation}
h(x,0) = h_0 + A \cos(2\pi x),
\label{h_cos}
\end{equation}
where $A=0.01$.

In order to compute the right-hand-side $f_j(h^n,t^n)$ of equation
\eqref{thin_film}, we use second-order centered finite differences
on a regular grid $x_j=j/N$, where $j\in[0,N]$ and $N=128$ is the number
of grid points~:
\begin{multline}
  f_j(h^n,t^n) = -h_j^3 \frac{h_{j-2}-4 h_{j-1} +
    6h_j - 4h_{j+1} + h_{j+2}}{\delta x^4}
	- 3 h_j^2 
	\frac{h_{j+1}-h_{j-1}}{2 \delta x} 
	\frac{-h_{j-2}+2h_{j-1}-2h_{j+1}+h_{j+2}}{2\delta x^3}\\
	- \frac{1}{h_j} \frac{h_{j+1}-2h_j+h_{j-1}}{\delta x^2}
	+ \frac{1}{h_j^2}\left( \frac{h_{j+1}-h_{j-1}}{2 \delta x} \right)^2.
	\label{thin_film_dis}
	\end{multline}
Using the Fourier transform of \eqref{thin_film_dis} in
\eqref{discrete_F_solve}, we obtain $\hat{h}^{n+1}_k$, from which 
the new points $h_j^{n+1}$ are obtained from the inverse Fourier transform.
The Richardson scheme \eqref{R}, based on each grid point, then leads 
to a second-order accurate result for $h_j^{n+1}$. 

\subsection{Stabilization}

Before proceeding to the automatic stabilization of the numerical scheme,
we adopt a von Neumann stability analysis, in order to predict the
theoretical value of $\lambda(k)$ for the scheme to be stable. 
For this purpose, we only need to consider the fourth-order derivative
in equation \eqref{thin_film}, which is the stiff term to be stabilized.
Using a ``frozen coefficient'' hypothesis, we look for perturbations 
to the mean profile $\overline{h}$ in the form of a single Fourier mode:
\beq
h_j^{n} = 
\overline{h}(j/N,n\delta t) + \delta \hat{h}_k^n
=
\overline{h}(j/N,n\delta t) + \xi^n e^{i k j \delta \alpha},
\label{sFm}
\eeq
where $\xi(\delta t,k)$ is the amplification factor \cite{NR07}, $\bar{h}$
is assumed constant over the time step,
$\delta \alpha = 2\pi \delta x = 2\pi / N$, and $k\in[0,N-1]$. 
Inserting this expression into \eqref{thin_film_dis},
retaining only $-h^3 h_{xxxx}$ in the original equation, 
the linearization  \eqref{linearize} of the right-hand-side of
\eqref{thin_film} gives~:
\begin{equation}
e(k) = 
2\frac{\overline{h}^3}{\delta x^4} 
\left(
\cos(2 k \delta \alpha)
-4 \cos(k \delta \alpha)
+3
\right).
\label{growth_TF}
\end{equation}

The modified numerical scheme \eqref{discrete_F_solve} will be stable
as long as $\lambda(k)$ meets the stability criterion \eqref{2ndstab} with
$e(k)$ given by \eqref{growth_TF}. Instead of this $\lambda(k)$,
for simplicity we initialize $\lambda(k)$ with an expression of the form
$\lambda_0 k^4$, where $\lambda_0$ must be chosen as the maximum of
$2e(k)/3 k^4$ over all $k$, in order to satisfy the stability criterion
everywhere. It is easily seen that this maximum is attained in the limit
$k\rightarrow 0$, from which we obtain 
\begin{equation}
\lambda_0 \ge \frac{32}{3}\pi^4 \overline{h}^3.
\label{small_k}
\end{equation}

Note that \eqref{small_k} depends on the local height $\overline{h}$,
which varies slightly for the initial condition; we choose it to
correspond to the maximum of the surface elevation, for which the
stability requirement is most stringent. In addition, the
condition \eqref{def_ke} defines a threshold value $k_e$,
below which the explicit time step is stable':
\begin{equation}
	\frac{2}{\delta t} = e(k) \simeq \overline{h}^3 (2\pi k_e)^4
	\Rightarrow 
	k_e \simeq \frac{1}{2\pi} \left( \frac{2}{ \delta t \;
          \overline{h}^3} \right)^{1/4} \simeq 4.25,
	\label{ke_thin}
\end{equation}
where we have used $\overline{h}=h_0=0.34$ and $\delta t=10^{-4}$.

\begin{figure}[htbp]
\begin{center}
\includegraphics{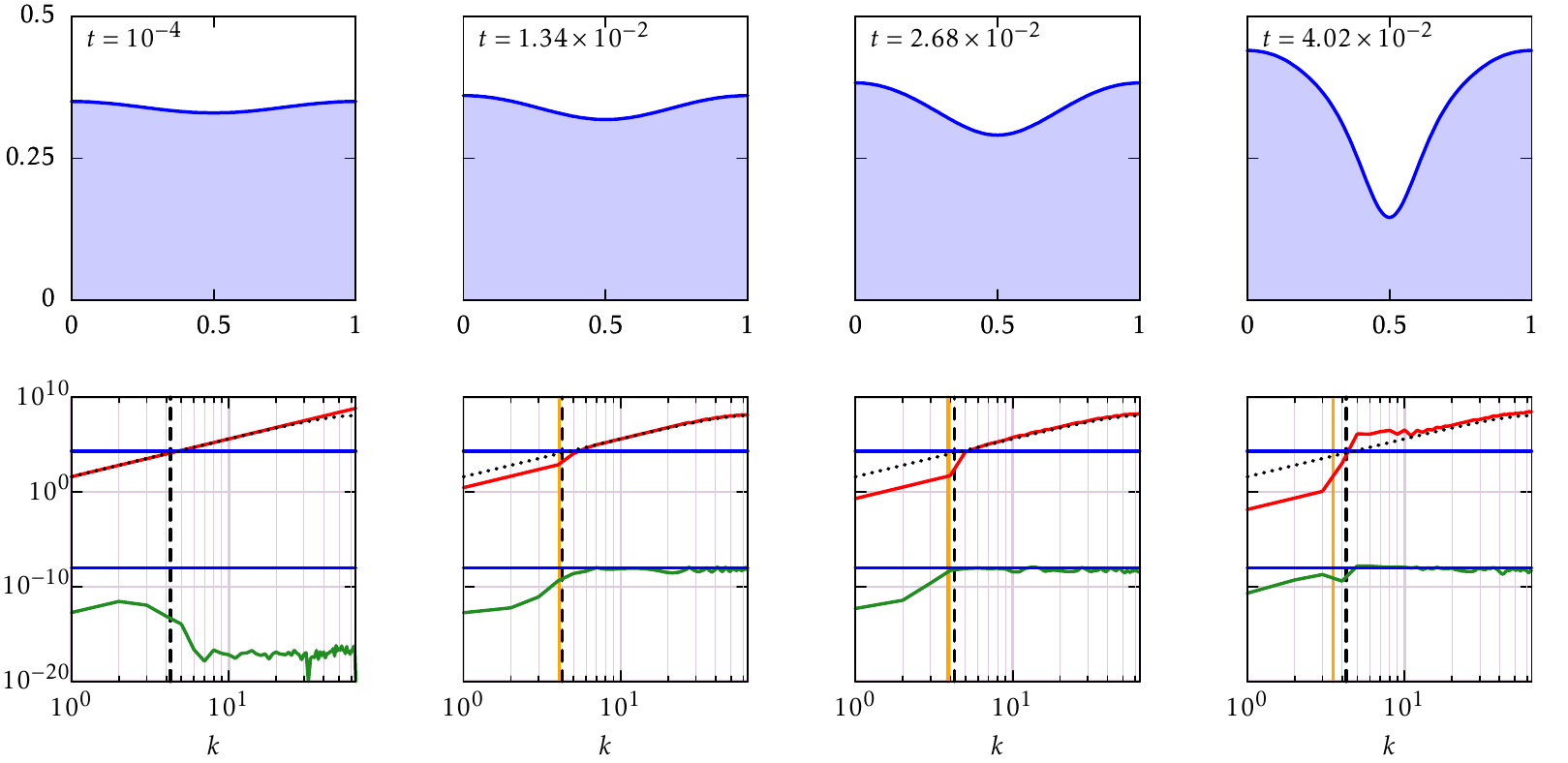}
\caption{
A simulation of the thin film equation with van der Waals forces
\eqref{thin_film}, with the interface 
shown on the top row. On the lower row, the corresponding spectrum of 
$\epsilon(k)$ defined by equation \eqref{epsilon} (green), as well as
$\lambda(k)$ (red). The dotted line is the stability
limit $\lambda_c(k) = 2e(k)/3$, with $e(k)$ given by \eqref{growth_TF}
and $\overline{h} \approx h_0$, the top horizontal blue line the explicit
stability boundary $2/\delta t$. The vertical dashed line is $k=k_e$
given by equation \eqref{ke_thin}, with $\overline{h} = h_0$, whereas the vertical orange line 
is $k=k_e$ computed for $\overline{h} = h_{max}$.
    }
    \label{evolution_thin_film}
\end{center}
\end{figure}

Using the procedure described in Sec.~\ref{sec:dynamic}, for each
Fourier mode $\epsilon(k)$ we adjust $\lambda(k)$ depending on whether
it is larger or smaller than an upper bound $\epsilon_u=10^{-8}$. 
The value of this bound is subject to some experimentation to make
the scheme work, but can be varied by several orders of magnitude without
affecting the functioning of the scheme. However, $\epsilon_u$ must be
chosen in accordance with the typical value of $\epsilon(k)$.
Indeed, in the current case, the absolute error \eqref{error} is
$\mathcal{O}(\delta t^2)$, and the noise measure \eqref{epsilon}
introduces a spatial error $\mathcal{O}(\delta x^2)$, leading to a global
absolute error $\mathcal{O}(\delta x^2 \delta t^2) \approx 6\times10^{-13}$,
consistent with the green curve in the first panel of
figure \ref{evolution_thin_film}. As a consequence $\epsilon_u$ has to be
larger than this value for the adaptive procedure to work.

Figure~\ref{evolution_thin_film} shows our adaptive scheme at work, as the
interface (shown on the top row) deforms; the noise measure $\epsilon(k)$ is
defined as in \eqref{epsilon} (See supplementary movie TF\_movie.mpeg). We initialized $\lambda(k)$ to the asymptotic
power-law $\lambda_0 k^4$, which is seen as the red line in the lower
panel of the first row (which shows the system after the first time step).
The true stability boundary, based on the full expression \eqref{growth_TF},
is shown as the dotted line; for small $k$, it is slightly lower than the
more stringent power law approximation chosen as the initial condition. 
The noise measure $\epsilon(k)$ after the first time step is very small
as expected.

As seen in the second panel, after some time $\lambda(k)$ has converged
onto the theoretical stability limit \eqref{2ndstab} for small $k$,
with $e(k)$ given by \eqref{growth_TF}. Here we have assumed
$\overline{h} \approx h_0$ for the initial stages of the dynamics,
an approximation that will no longer be valid near the end of the
computation (fourth row of Figure~\ref{evolution_thin_film}), where
$h$ varies considerably in space, and \eqref{sFm} can only be applied
locally. However, adjustment of $\lambda(k)$ toward the stability
boundary only occurs for $k>k_e$, since there is no numerical
instability below $k=k_e$. As a result, for $k<k_e$ the stabilizing
spectrum $\lambda(k)$ is reduced at every time step, and in the second
panel has already fallen by orders of magnitude below
the stability limit of the EIN scheme. 

The only source for concern is seen in the 4th panel, when the film
thickness has become very non-uniform. In that case there is a
region just above $k_e$ where $\lambda(k)$ is quite elevated relative to
the theoretical limit (but which is based on the assumption of a uniform
thickness), as well as noisy. 
This comes from the fact that the explicit stability boundary $k=k_e$
given by \eqref{growth_TF} moves to the left when values of $h$ in space
become significantly higher than $h_0$. The vertical orange line in
figure \ref{evolution_thin_film} corresponds to $k=k_e$ computed using the
maximum height of the interface $h_{max}$ instead of $h_0$. As long as
this value does not reach the next smaller integer value of $k$,
$\lambda(k)$ is progressively decreased on the left of this boundary. 
As soon as $k_e$ crosses an integer value $k'$, exponential growth is
observed for $k'$ and $\lambda(k')$ cannot be adapted quickly enough,
resulting in increasing values of $\lambda$ for the adjacent values
of $k$. This problem could probably be solved by changing the way
$\lambda(k)$ is adapted at each timestep, by testing the stability of
each mode and rejecting the timestep in case of instability.


\section{Example: 2D Kuramoto--Sivashinsky equation}
\label{sec_kur}

\subsection{Equation of motion}

\begin{figure}[htbp]
\begin{center}
\includegraphics{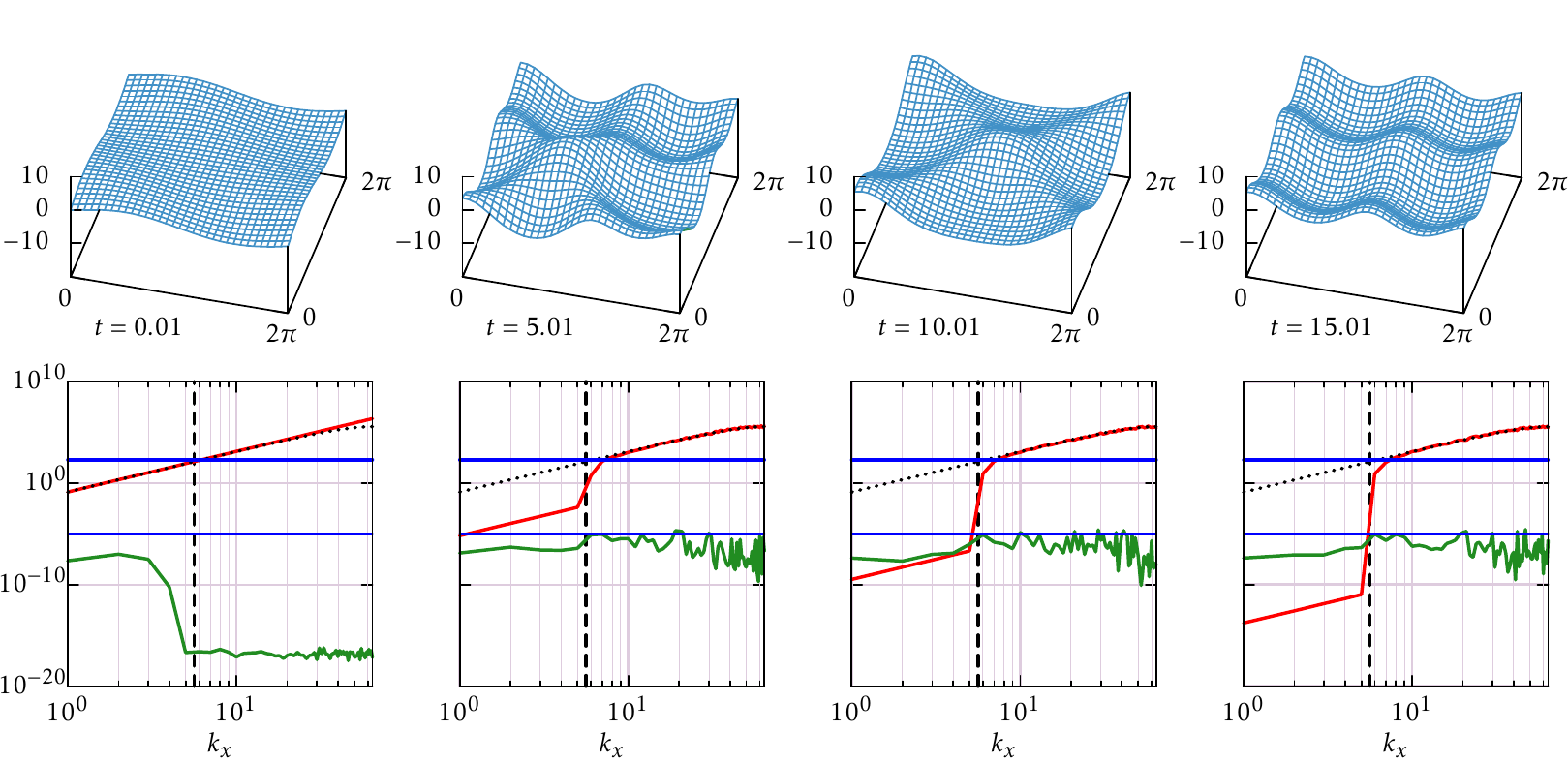}
\caption{
A simulation of the Kuramoto--Sivashinsky equation \eqref{KS} for
$\nu=0.2$, with the interface shown on the top row. 
$N_x=N_y=128$ points are used to discretize \eqref{KS} in both
directions, and the time step is $\delta t = 0.01$.
On the lower row, the corresponding spectrum of 
$\epsilon(k_x,0)$ defined by equation \eqref{epsilon} (green), as well as
$\lambda(k_x,0)$ (red). The dotted line is the stability
limit $\lambda_c(k_x,0) = 2e(k_x,0)/3$ (cf. \eqref{growth_KS});
the top horizontal blue line the explicit stability boundary $2/\delta t$. 
The vertical dashed line is $k=k_e$ given by $2/\delta t = e_s(k_x,0)$ (cf. \eqref{smallk_KS});
$\lambda(k_x,k_y)$ is initialized using $2e_s(k_x,k_y)/3$.
	Half of the spectrum in $k_x$ is shown, since it is symmetric around $N_x/2$ (the other half corresponding to negative wave-numbers $k_x \in [-N_x/2 + 1:-1]$).
    }
    \label{evolution_KS}
\end{center}
\end{figure}

To demonstrate that our method works in higher dimensions, we consider 
the example of the 2D Kuramoto--Sivashinsky equation \cite{cross93},
which is known to exhibit spatio-temporal chaos~\cite{cross93,kalogirou2015}~:
\begin{equation}
  \frac{\partial u}{\partial t} =
  -\mathcal{N}(u) - \Delta u - \nu \Delta^2 u .
\label{KS}
\end{equation}
The single Laplacian on the right has a minus sign in front of it,
leading to instability on the smallest scale; this is stabilized by the
last term, which is of forth order, making the problem very stiff. 
Nonlinearity is introduced through the term 
\begin{equation}
  \mathcal{N}(u) = \frac{1}{2} \left( | \nabla u |^2
  - \frac{1}{4 \pi^2} \int_0^{2 \pi} \int_0^{2 \pi} | \nabla u |^2
  dx dy \right); 
\end{equation}
following \cite{kalogirou2015}, the spatially constant integral term has
been introduced for convenience only, to make sure that $u$ always has
zero spatial mean.
The variable $u(x,y,t)$ is defined on a two-dimensional square 
domain, which we can rescale to ensure that $(x,y) \in [0,2 \pi]$.

Equation \eqref{KS} is discretized on a regular grid using centered
finite differences in order to find $f_{j_x,j_y}(u^n,t^n)$, whose two-dimensional
Fourier transform is $\hat{f}_{k_x,k_y}(u^n,t^n)$. We can then use the
modified time step \eqref{discrete_F_solve} to find
$\hat{u}^{n+1}_{k_x,k_y}$ and thus $u_{j_x,j_y}^{n+1}$. As usual, the scheme is
then turned into a second order method using Richardson extrapolation
\eqref{R}. 
In \cite{kalogirou2015}, \eqref{KS} is treated implicitly using a
Fourier pseudospectral method. The purpose of our treatment is to demonstrate
the effectiveness of our general scheme, which does not pay attention to
the specifics of the operator, in spatial dimensions greater than one.

\subsection{Stabilization}
In order to study the numerical stability, we only consider the
bi-Laplacian, which is the stiffest term to be stabilized. 
Inserting the single Fourier mode 
$$
u_{j_x,j_y}^n = \xi^n e^{i(k_x j_x \delta x + k_y j_y \delta y)}
$$
into the discretized version of equation \eqref{KS}, and retaining
only the bi-Laplacian term, we obtain~:
\begin{eqnarray}
	\nonumber e(k_x,k_y) & = &  \nu \left\{ \frac{2 \cos 2 k_x \delta x - 8 \cos k_x \delta x + 6}{\delta x^4} + \frac{2 \cos 2 k_y \delta y - 8 \cos k_y \delta y + 6}{\delta y^4} \right. \\
&& +  \left. 2 \frac{2 \cos (k_x \delta x + k_y \delta y) -
  4 \cos k_y \delta y + 2 \cos (k_x \delta x - k_y \delta y) -
	4 \cos k_x \delta x + 4}{\delta x^2 \delta y^2} \right\} . 
\label{growth_KS} 
\end{eqnarray}

The modified numerical scheme \eqref{discrete_F_solve} will be stable
as long as $\lambda(k_x,k_y)$ meets the stability criterion \eqref{2ndstab}~: 
$\lambda(k_x,k_y) > 2e(k_x,k_y)/3$, and the stability of the explicit
scheme is given by \eqref{def_ke}: $e(k_x,k_y) < 2/\delta t$. 
For the explicit stability boundary, one can take the small-$k$ limit 
of \eqref{growth_KS}: 
\beq
e_s(k_x,k_y) \sim  \nu \left( k_x^2 + k_y^2 \right)^2. 
\label{smallk_KS}
\eeq
As in the previous example, this approximation overpredicts the true value 
of $e(k_x,k_y)$ for large wavenumbers. For simplicity, we will use this
conservative approximation to set the initial value of $\lambda(k_x,k_y)$.
The adjustment of $\lambda$ is based on an upper bound $\epsilon_u=10^{-5}$
for $\epsilon$. 
This value is significantly higher than in the previous example, 
since the typical absolute error in the present case is 
$\mathcal{O}(\delta x^2 \delta t^2)\approx 2\times 10^{-7}$, as seen on the green curve in the first panel of figure \ref{evolution_KS}.

In the definition \eqref{error} of $\epsilon(k_x,k_y)$, we now have to
interpolate the error estimator $E(j_x,j_y)$ on a two-dimensional grid;
we find $\bar{E}(j_x,j_y)$ using a bilinear interpolation of the error
estimator~:
\beq
\bar{E}(j_x,j_y) = \frac{1}{4} \left[ E(j_x,j_y-1)+E(j_x+1,j_y)+E(j_x,j_y+1)+E(j_x-1,j_y) \right]
\eeq

\begin{figure}[htbp]
\begin{center}
\includegraphics{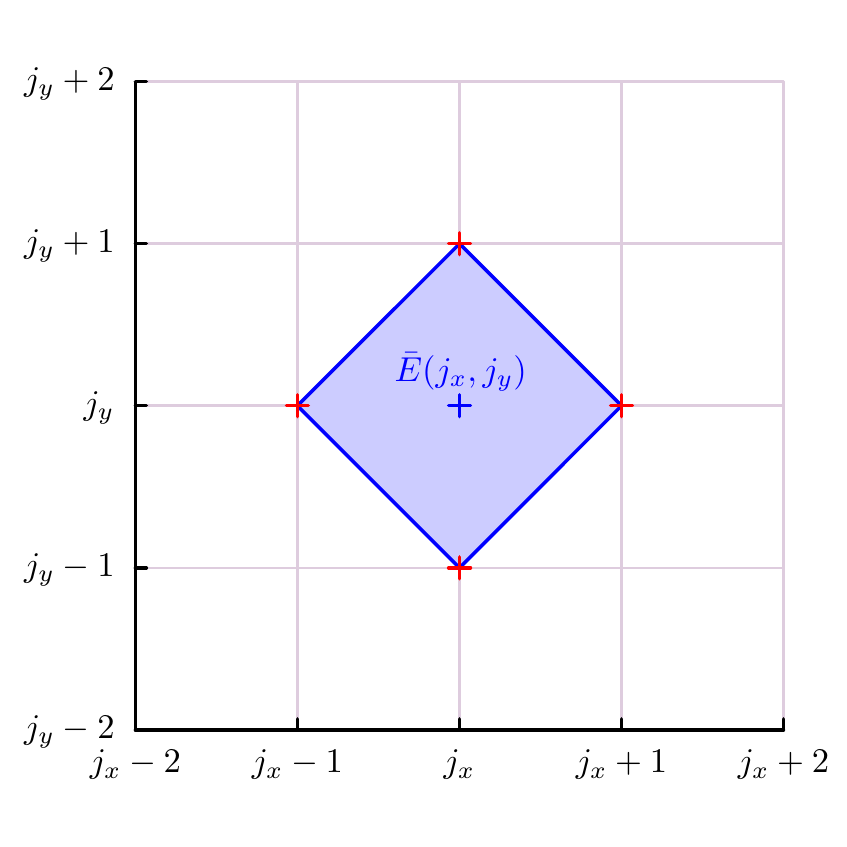}
\caption{
	Bilinear interpolation $\bar{E}(j_x,j_y)$ using the four neighboring points in red. 
  }
    \label{bilint}
\end{center}
\end{figure}

Our first attempt was to use $E(j_x+1,j_y-1)$, $E(j_x+1,j_y+1)$, $E(j_x-1,j_y+1)$ and $E(j_x-1,j_y-1)$ to interpolate $E$ in $(j_x,j_y)$, 
but it turned out that the interpolation values for two adjacent nodes where decoupled, making the process of estimating the error unstable.
This issue is addressed by using the four neighbors $E(j_x,j_y-1)$, $E(j_x+1,j_y)$, $E(j_x,j_y+1)$, and $E(j_x-1,j_y)$, as seen in figure \ref{bilint}.

Fig.~\ref{evolution_KS} shows a computation of the Kuramoto--Sivashinsky
equation \eqref{KS}, for $\nu=0.2$, in the chaotic regime (See supplementary movie KS\_movie.mpeg). Thus
the interface $u(x,y)$ deforms in an irregular, unpredictable fashion
on many scales. On the lower row we show the corresponding spectrum of
$\epsilon(k_x,0)$ (green), as well as $\lambda(k_x,0)$ (red). For plotting
purpose, we chose to show only a slice of the spectrum ($k_y=0$), but
the adaption procedure works for the whole 2D spectrum.

We have initialized $\lambda(k_x,k_y)$ to $2e_s(k_x,k_y)/3$ (see
\eqref{smallk_KS}). As in the previous example, the initial condition
for $\lambda$ (red line) is slightly above the theoretical stability
limit (dotted line) for $k$ values corresponding to small scales. 
This is still true after the first time step, while $\epsilon(k_x,0)$
is very small as expected. As seen in the second panel, $\lambda(k_x,0)$
has converged onto the theoretical stability limit
$\lambda_c(k_x,0) = 2e(k_x,0)/3$, with $e(k_x,0)$ given by
\eqref{growth_KS}, since the initial condition overpredicts the stability
boundary.

However, convergence only occurs for $k>k_e$, since below $k=k_e$ no
numerical instability occurs. As a
result, for $k<k_e$ the stabilizing spectrum $\lambda(k_x,0)$ is reduced
at every time step, and has already fallen by orders of magnitude below
the stability limit of the EIN scheme. Correspondingly, by adjusting
$\lambda(k_x,0)$ the error $\epsilon(k_x,0)$ is kept close to the threshold
$\epsilon_u=10^{-5}$ for $k>k_e$.


\section{Example: Hele--Shaw flow}
\label{sec_HS}

\subsection{Equations of motion}

As an example of a non-local, but stiff operator, 
we consider an interface in a vertical Hele-Shaw cell, separating two
viscous fluids with the same dynamic viscosity, with the heavier fluid
on top \cite{HLS94}. As heavy fluid falls, small perturbations on the
interface grow exponentially: this is known as the Rayleigh-Taylor
instability \cite{Draz81}. However, surface tension assures regularity
on small scales. For simplicity, we assume the flow to be periodic in the
horizontal direction. We briefly recall the dynamics of the interface
here; for more details, see \cite{HLS94,DE14}. 

The interface is discretized using marker points labeled with $\alpha$, 
which represents the motion of a fluid particle. They are advected
according to~:
\beq
\frac{\partial{\bf X}(\alpha)}{\partial t} = U {\bf n} + T {\bf s}.
\label{advection}
\eeq
Here ${\bf X}(\alpha)=(x,y)$ is the position vector,
${\bf n}=(-y_\alpha/s_\alpha,x_\alpha/s_\alpha)$ and
${\bf s}=(x_\alpha/s_\alpha,y_\alpha/s_\alpha)$ are the normal and tangential
unit vectors, respectively, and $s_\alpha=(x_\alpha^2+y_\alpha^2)^{1/2}$.
Hence $U=(u,v)\cdot {\bf n}$ and $T=(u,v)\cdot {\bf s}$ are the normal
and tangential velocities, respectively. The tangential velocity does not
affect the motion, but is chosen so as to maintain a reasonably uniform
distribution of points \cite{HLS94,DE14}.
If $z(\alpha,t)=x+iy$ is the complex position of the interface, which
is assumed periodic with period $1$ ($z(\alpha+2\pi)=z(\alpha)+1$), the
complex velocity becomes~:
\beq
u(\alpha)-i v(\alpha)=\frac{1}{2i} 
PV\int_{0}^{2\pi}{\gamma(\alpha',t)\cot
\left[ \pi(z(\alpha,t)-z(\alpha',t))\right] d\alpha'},
\label{BR}
\eeq
where $\gamma$ is the vortex sheet strength. For two fluids of equal 
viscosities \cite{MB02},
\beq
\gamma = S \kappa_{\alpha} - R y_\alpha,
\label{gamma}
\eeq
where $\kappa$ is the mean curvature of the interface~:
\beq
\kappa(\alpha)=\frac{x_\alpha y_{\alpha \alpha}-y_\alpha x_{\alpha \alpha}}
{s_\alpha^3}, \qquad \mathrm{recalling \,\,\,\, that} \qquad
s_\alpha=(x_\alpha^2+y_\alpha^2)^{1/2}.
\label{kappa}
\eeq
Here $S$ is the non-dimensional surface tension coefficient and $R$ is
the non-dimensional gravity force, chosen to be $0.1$ and $50$,
respectively, in the following example. To compute the complex Lagrangian
velocity of the interface \eqref{BR}, we use the spectrally accurate
alternate point discretization \cite{S92}~:
\beq
u_j - i v_j \simeq -\frac{2\pi i}{N} \sum_{{l=0}\atop{j+l \; odd}}^{N-1}
{\gamma_l \cot\left[\pi(z_j-z_l)\right]}.
\label{comp_vel}
\eeq
Derivatives $\kappa_\alpha$ and $y_\alpha$ are computed at each time step using 
second-order centered finite differences, and $\alpha$ is defined by 
$\alpha(j) = 2 \pi j/N$, where $j \in [0,N]$ and $N=1024$ is the number of 
points describing the periodic surface. Note that the numerical effort
of evaluating (\ref{comp_vel}) requires $\mathcal{O}(N^2)$ operations,
and thus will be the limiting factor of our algorithm. 

\subsection{Stabilization}

Although the character of the non-local operator in this example is
very different from previous equations, our numerical stabilization
works in a fashion that is remarkably similar. The modified scheme
\eqref{discrete_F_solve} now becomes
\beq
\hat{x}^{n+1}_k = \hat{x}^{n}_k + 
\frac{\hat{u}_k^n}{
\delta t^{-1} + \lambda(k)}, \quad
\hat{y}^{n+1}_k = \hat{y}^{n}_k + 
\frac{\hat{v}_k^n}{
\delta t^{-1} + \lambda(k)},
\label{xy_equ}
\eeq
where $\hat{u}_k^n$ and $\hat{v}_k^n$ are calculated from the Fourier
transform of (\ref{comp_vel}). The new grid points $x_j^{n+1},y_j^{n+1}$
are obtained from the inverse Fourier transform of
$\hat{x}^{n+1}_j,\hat{y}^{n+1}_j$, and for each
component \eqref{xy_equ} is turned into a second-order scheme using 
\eqref{R}.

In \cite{DE14}, we performed a linear analysis \eqref{linearize} of
the discrete modes of \eqref{comp_vel} about a flat interface. We found that
\beq
e(k) = 
\frac{S N^3}{L^3}\left(1 - \cos\frac{2\pi k}{N}\right)\sin\frac{2\pi k}{N}
\equiv \tilde{e}(x) = \left(1 - \cos x\right)\sin x, \quad
x = \frac{2\pi k}{N},
\label{growth_HS}
\eeq
where $L$ is the length of the interface, and $N$ the number of gridpoints.
As before, we use the long-wavelength approximation to $e(k)$:
\beq
e(k) \approx \frac{S}{2L^3} (2\pi k)^3 \equiv e_s(x) = \frac{x^3}{2},
\label{a_power}
\eeq
to find the explicit stability boundary \eqref{def_ke} as 
\beq
k_e \approx \frac{L}{2\pi} \left(\frac{4}{S\delta t}\right)^{1/3}.
\label{k_expl}
\eeq

\begin{figure}[htbp]
\begin{center}
\includegraphics[width=\linewidth]{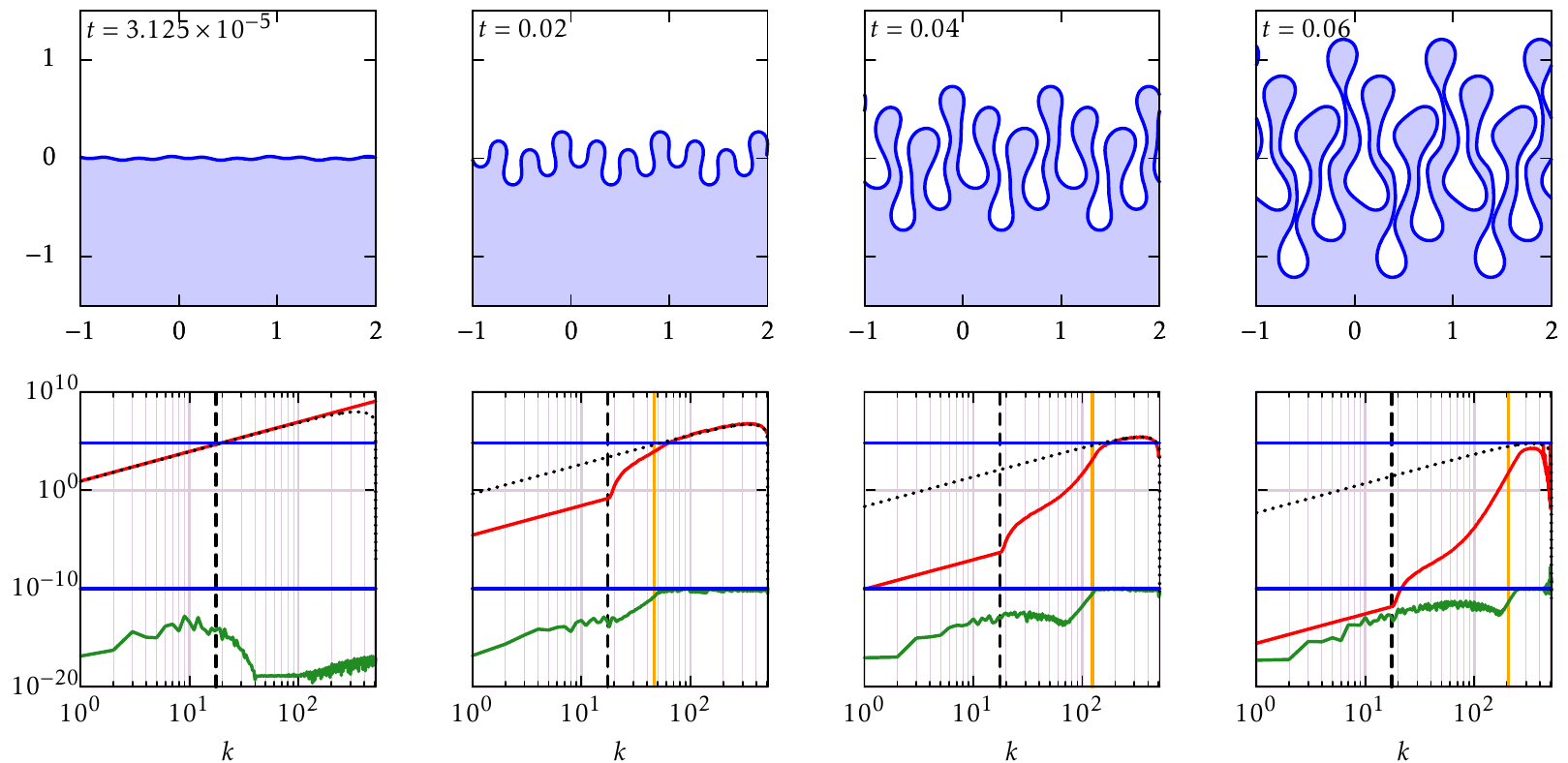}
\caption{
A simulation of the Hele-Shaw problem \eqref{advection}, \eqref{BR}
shown on the top row. On the lower row, the corresponding spectrum of 
$\epsilon(k)$ defined by equation \eqref{errorHS} (green), as well as
$\lambda(k)$ (red). The dotted line is the stability
limit $\lambda_c(k) = 2e(k)/3$, with $e(k)$ given by \eqref{growth_HS},
the top horizontal blue line the explicit stability boundary $2/\delta t$.
The vertical orange line is $k=k_e$ given by equation \eqref{k_expl}, and 
the vertical dashed line its initial position.
    }
\label{HS_stable}
\end{center}
\end{figure}

Using the same approximation \eqref{a_power}, which overpredicts
the critical value $\lambda(k)$, one finds 
\beq
\lambda_c(k) = \frac{S}{3} \left(\frac{2\pi k}{L}\right)^3
\label{theor_stab}
\eeq
as a sufficient condition for stability. In \cite{DE14}, we used
a fixed spectrum $\lambda(k)$, slightly larger than \eqref{theor_stab},
to stabilize the Hele-Shaw dynamics. 

We now use the same procedure as before, with the upper bound
$\epsilon_u=10^{-10}$ for $\epsilon(k)$. 
Figure~\ref{HS_stable} shows our adaptive scheme at work, as the
interface (shown on the top row) deforms, and the length $L$ of the
interface increases (See supplementary movie HS\_movie.mpeg). The error $\epsilon(k)$ is defined as the
maximum of \eqref{epsilon} over the two components:
\beq
\epsilon(k) = MAX(\hat{E}^x_k - \hat{\bar{E}}^x_k,
\hat{E}^y_k - \hat{\bar{E}}^y_k).
\label{errorHS}
\eeq

Initializing $\lambda(k)$ to the approximation \eqref{theor_stab}
(red line), we observe the same convergence toward the theoretical
stability boundary as before (dotted line). A new feature is that
on account of the length $L$ of the boundary increasing in time,
the explicit stability boundary $k_e$ (vertical dashed and orange lines) increases
in time, and the theoretical stability boundary for $\lambda(k)$
(dotted line) comes down. As seen in the second panel of
Fig.~\ref{HS_stable}, our adaptive scheme for $k > k_e$ has converged
toward the theoretical prediction, and then continues to trace
it as $L$ increases. The region where no stabilization is required
increases as well, and $\lambda(k)$ decreases to very low values on an
increasingly large domain. 
The results described above are not changed significantly as
$\epsilon_{u}$ is varied over several orders of magnitude up or down
from $10^{-10}$, but of course the value must be significantly over the
rounding error, and below the expected truncation error. 

\begin{figure}[htbp]
\begin{center}
\includegraphics[width=\linewidth]{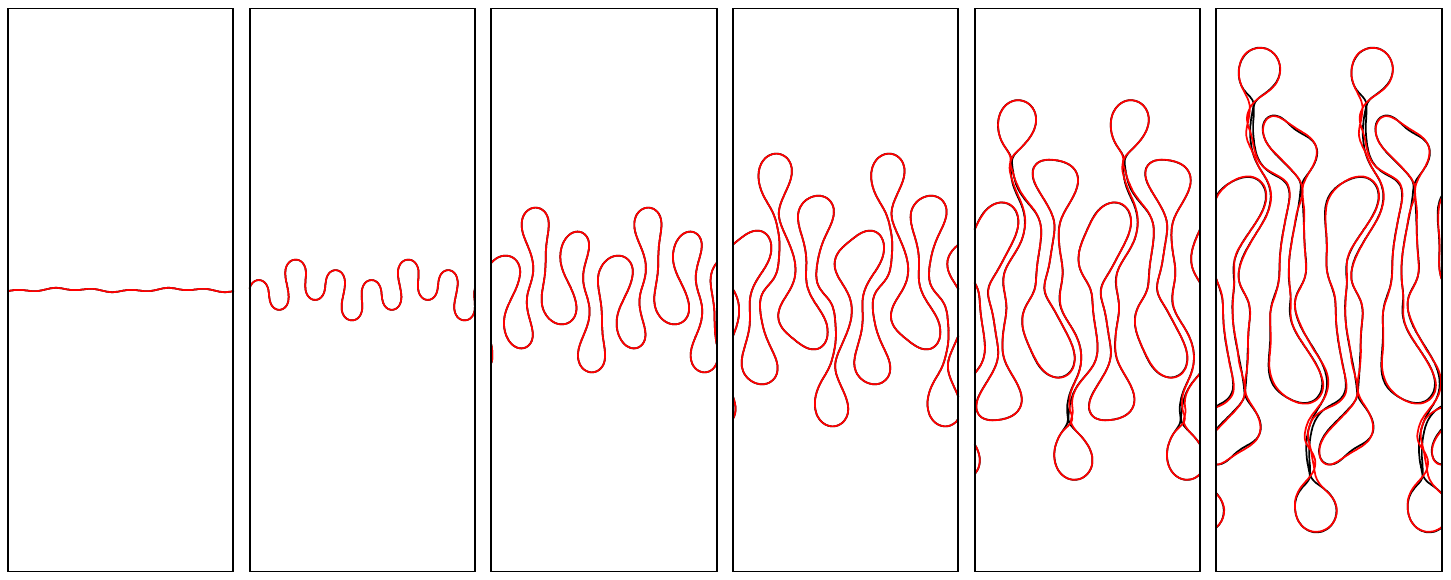}
\caption{A comparison of the interface as obtained from our current
adaptive scheme (red curves) and our earlier EIN scheme (black curves)
\cite{DE14}, which used the theoretical stability boundary \eqref{theor_stab}. 
 }
\label{HS_compare}
\end{center}
\end{figure}

In our earlier EIN scheme \cite{DE14}, we used $\lambda(k)$ based on the
simplified stability boundary \eqref{theor_stab} to stabilize
the Hele-Shaw interface motion shown in Fig. \ref{HS_stable}. However,
this overpredicts the necessary damping for large $k$. In addition,
for $k<k_e$, no damping is necessary, and our adaptive scheme
reflects that by decreasing $\lambda(k)$ more and more. As a result,
the damping in the adaptive scheme is significantly smaller than
in our previous EIN scheme.
In Fig.~\ref{HS_compare} we show a comparison of the numerical results
to those of the earlier scheme, and find very good agreement. The major
advance is of course that $\lambda(k)$ no longer needs to be prescribed,
but is found self-consistently as part of the algorithm which ensures
stability. Only in the last panel is there a significant discrepancy
between the two results. This occurs in places where two sides of the
interface have come in close proximity, comparable to the spacing
between grid points. But this means our evaluation of the velocity integral
is no longer sufficiently accurate to be reliable. 
                                   
\section{Outlook and conclusions}

We have demonstrated the feasibility of our method using three 
different model problems, highlighting different aspects of physical
problems containing a wide spectrum of time scales, making them stiff.
Clearly, there are many ways in which to extend and
improve the present approach. Firstly, we estimated the damping spectrum 
by analyzing the current solution in Fourier space, which is particularly
easy for the periodic domain considered by us. However, this may be
circumvented by periodically continuing a solution defined over a finite
domain only. In addition, one could formulate the entire method in real
space, as done in some cases described in \cite{DE14}. 

A second, more important issue is our assumption of the spectrum
$e(k)$ in \eqref{linearize} being real. This assumption is well founded,
since the ultimate physical damping process is dissipative, leading to
real eigenvalues. However, as demonstrated by the example of an inertial
vortex sheet considered in \cite{HLS94}, even problems lacking dissipation
can display significant stiffness. This case leads to a system of PDEs,
with pairs of complex eigenvalues $e(k)$ on the right-hand-size of
\eqref{linearize}, corresponding to traveling waves. In that case the
damping spectrum $\lambda(k)$ would also have to be complex to ensure
stability \cite{DE14}, a case we have not yet considered. 

Finally, a problem we still need to address is how to choose an
initial condition for the damping spectrum $\lambda(k)$. In the
present work we choose a power-law spectrum which can be inferred
from a simple analysis of the continuum version of the equations of
motion, which then adapts to an optimal spectrum. It would be ideal
if no input whatsoever was necessary, choosing for example
$\lambda(k)=0$ initially. At present, this is not possible, as the
quality of the numerical solution deteriorates before $\lambda(k)$
can adapt. We suspect that in order for such a scheme to be successful,
one needs to implement a variable time step, such that initial steps
during which $\lambda(k)$ is found are very small. 
                       
In conclusion, following our previous study on this subject,
we propose a new method to remove the stiffness of PDEs containing
non-linear stiff terms, {\it i.e.} high spatial derivatives embedded
into non-linear terms. This method allows for the self-consistent
estimation of a stabilizing term on the right-hand-side of the PDE, that
ensures absolute stability for the numerical scheme. Analyzing the
spectrum of the solution at each time step, we adapt automatically
the stabilizing term such that each unstable Fourier mode is damped optimally. 
The computational cost of this method is essentially the same as that
of the explicit method. 


\bibliographystyle{model1-num-names}
\bibliography{../all_ref.bib}

\end{document}